\documentclass{amsart}

\usepackage{amsfonts,amsmath,amssymb,mathrsfs,verbatim,cite}
\usepackage[enableskew]{youngtab}

\numberwithin{equation}{section}

\newtheorem{lemma}{Lemma}[section]
\newtheorem{theorem}{Theorem}[section]
\newtheorem{corollary}{Corollary}[section]

\newcommand{\abs}[1]{\lvert#1\rvert}

\newcommand{\F}{\mathscr F}
\newcommand{\G}{\mathscr G}
\newcommand{\Z}{\mathbb Z}
\newcommand{\la}{\lambda}
\newcommand{\lab}{\bar{\lambda}}
\newcommand{\mub}{\bar{\mu}}
\newcommand{\omegab}{\bar{\omega}}
\newcommand{\thetab}{\bar{\theta}}

\newcommand{\qbin}[2]{\genfrac{[}{]}{0pt}{}{#1}{#2}}
\newcommand{\pp}[2]{(#1|#2)}
\newcommand{\vul}{\bullet}

\begin{document}

\title[Hall--Littlewood functions]
{Hall--Littlewood functions and the A$_2$ Rogers--Ramanujan identities}

\author{S. Ole Warnaar}\thanks{Work supported by the Australian
Research Council}
\address{Department of Mathematics and Statistics,
The University of Melbourne, VIC 3010, Australia}
\email{warnaar@ms.unimelb.edu.au}



\begin{abstract}
We prove an identity for Hall--Littlewood symmetric functions labelled
by the Lie algebra A$_2$. Through specialization this yields a
simple proof of the A$_2$ Rogers--Ramanujan identities of Andrews, Schilling
and the author.
\end{abstract}

\maketitle

\section{Introduction}
The Rogers--Ramanujan identities, given by \cite{Rogers94}
\begin{subequations}\label{RR}
\begin{equation}\label{RR1}
1+\sum_{n=1}^{\infty}\frac{q^{n^2}}{(1-q)(1-q^2)\cdots (1-q^n)}
=\prod_{n=1}^{\infty}\frac{1}{(1-q^{5n-1})(1-q^{5n-4})}
\end{equation}
and
\begin{equation}
1+\sum_{n=1}^{\infty}\frac{q^{n(n+1)}}{(1-q)(1-q^2)\cdots (1-q^n)}
=\prod_{n=1}^{\infty}\frac{1}{(1-q^{5n-2})(1-q^{5n-3})},
\end{equation}
\end{subequations}
are two of the most famous $q$-series identities, with deep connections
with number theory, representation theory, statistical mechanics
and various other branches of mathematics.

Many different proofs of the Rogers--Ramanujan identities have been given
in the literature, some bijective, some representation theoretic,
but the vast majority basic hypergeometric.
In 1990, J.~Stembridge, building on work of I.~Macdonald, 
found a proof of the Rogers--Ramanujan identities quite unlike any of the
previously known proofs. In particular he discovered that Rogers--Ramanujan-type
identities may be obtained by appropriately specializing identities
for Hall--Littlewood polynomials. 
The Hall--Littlewood polynomials and, more generally, Hall--Littlewood functions
are an important class of symmetric functions,
generalizing the well-known Schur functions. Stembridge's Hall--Littlewood
approach to Rogers--Ramanujan identities has been further generalized in
recent work by Fulman \cite{Fulman00},
Ishikawa \textit{et al.} \cite{IJZ04} and Jouhet and Zeng \cite{JZ04}.

Several years ago Andrews, Schilling and the present author generalized
the two Rogers--Ramanujan identities to three identities labelled by the
Lie algebra A$_2$ \cite{ASW99}. 
The simplest of these, which takes the
place of \eqref{RR1} when A$_1$ is replaced by A$_2$ reads
\begin{align}\label{RRA2}
(q;q)_{\infty}&\sum_{n_1,n_2=0}^{\infty}\frac{q^{n_1^2-n_1 n_2+n_2^2}}
{(q;q)_{n_1}(q;q)_{n_2}(q;q)_{n_1+n_2}} \\
&=\sum_{n_1,n_2=0}^{\infty}\frac{q^{n_1^2-n_1 n_2+n_2^2}}
{(q;q)_{n_1}}\qbin{2n_1}{n_2} \notag \\
&=\prod_{n=1}^{\infty}
\frac{1}{(1-q^{7n-1})^2(1-q^{7n-3})(1-q^{7n-4})(1-q^{7n-6})^2},
\notag
\end{align}
where
$(q;q)_0=1$ and $(q;q)_n=\prod_{i=1}^n (1-q^i)$ is a $q$-shifted factorial, and 
\begin{equation*}
\qbin{n}{m}=\qbin{n}{m}_q=
\begin{cases}
\displaystyle
\frac{(q^{n-m+1};q)_m}{(q;q)_m} & \text{for $m\geq 0$,} \\[3mm]
0 & \text{otherwise}
\end{cases}
\end{equation*}
is a $q$-binomial coefficient.
The equivalence of the two expressions 
on the left of \eqref{RRA2} follows from a straightforward application
of Jackson's terminating $_2\phi_1$ transformation
\cite[Equation (III.7)]{GR90}, see \cite{ASW99}.

The A$_2$ characteristics of \eqref{RRA2} are (i) the exponent of $q$ of the
two summands, which may alternatively be put as
$\frac{1}{2}\sum_{i,j=1}^2 C_{ij}n_in_j$ with
$C=((2,-1),(-1,2))$ the A$_2$ Cartan matrix, and (ii)
the infinite product on the right, which can be identified with a 
branching function of the coset pair 
$(\textup{A}_2^{(1)}\oplus\textup{A}_2^{(1)},\textup{A}_2^{(1)})$ 
at levels $-9/4$, $1$ and $-5/4$, see \cite{ASW99}.

An important question is whether \eqref{RRA2} and its companions can again
be understood in terms of Hall--Littlewood functions.
This question is especially relevant since the A$_n$ analogues of the
Rogers--Ramanujan identities have so far remained elusive, and
an understanding of \eqref{RRA2} in the context of symmetric
functions might provide further insight into the structure of the
full A$_n$ generalization of \eqref{RR}.

In this paper we will show that the theory of
Hall--Littlewood functions may indeed be applied to 
yield a proof of \eqref{RRA2}. In particular we will
prove the following A$_2$-type identity for Hall--Littlewood functions.
\begin{theorem}\label{mainthm}
Let $x=(x_1,x_2,\dots)$, $y=(y_1,y_2,\dots)$ and let
$P_{\la}(x;q)$ and $P_{\mu}(y;q)$ be Hall--Littlewood functions
indexed by the partitions $\la$ and $\mu$. Then
\begin{multline}\label{Pid}
\sum_{\la,\mu} q^{n(\la)+n(\mu)-\pp{\la'}{\mu'}} P_{\la}(x;q)P_{\mu}(y;q) \\
=\prod_{i\geq 1}\frac{1}{(1-x_i)(1-y_i)}
\prod_{i,j\geq 1}\frac{1-x_iy_j}{1-q^{-1}x_iy_j}.
\end{multline}
\end{theorem}
In the above $\la'$ and $\mu'$ are the conjugates of $\la$ and $\mu$,
$\pp{\la}{\mu}=\sum_{i\geq 1}\la_i\mu_i$, 
and $n(\la)=\sum_{i\geq 1}(i-1)\la_i$.

For $q=1$ the Hall--Littlewood function $P_{\la}(x;q)$ reduces to the monomial
symmetric function $m_{\la}(x)$, and the identity \eqref{Pid} factorizes
into a product of the well-known
\begin{equation*}
\sum_{\la}m_{\la}(x)=\prod_{i\geq 1}\frac{1}{1-x_i}.
\end{equation*}

An appropriate specialization of Theorem~\ref{mainthm} leads to
a $q$-series identity of \cite{ASW99} which is the key-ingredient in
proving \eqref{RRA2}. In fact, the steps leading from \eqref{Pid} to \eqref{RRA2}
suggests that what is needed for the A$_n$ version of the Rogers--Ramanujan
identities is an identity for the more general sum
\begin{multline}\label{Ansum}
\sum_{\la^{(1)},\dots,\la^{(n)}} 
\prod_{i=1}^n q^{n(\la^{(i)})-\pp{{\la^{(i)}}'}{{\la^{(i+1)}}'}}
P_{\la^{(i)}}(x^{(i)};q),
\end{multline}
where $\la^{(1)},\dots,\la^{(n+1)}$ are partitions with
$\la^{(n+1)}=0$ the empty or zero partition,
and $x^{(i)}=(x_1^{(i)},x_2^{(i)},\dots)$.
What makes this sum difficult to handle is that
no factorized right-hand side exists for $n>2$.

\medskip

In the next section we give the necessary background material on 
Hall--Littlewood functions.
In Section~\ref{secC} some immediate consequences of Theorem~\ref{mainthm}
are derived, the most interesting one being the new
$q$-series identity claimed in Corollary~\ref{Cor4}.
Section~\ref{secP} contains a proof of Theorem~\ref{mainthm} and
Section~\ref{secRR} contains a proof of the A$_2$ Rogers--Ramanujan identities
\eqref{RRA2} based on Corollary~\ref{Cor4}.
Finally, in Section~\ref{secF} we present some open problems
related to the results of this paper.

\section{Hall-Littlewood functions}
We review some basic facts from the theory of Hall-Littlewood functions.
For more details the reader may wish to consult Chapter III of Macdonald's
book on symmetric functions \cite{Macdonald95}.

Let $\lambda=(\lambda_1,\lambda_2,\dots)$ be a partition, 
i.e., $\lambda_1\geq \lambda_2\geq \dots$ with finitely many $\lambda_i$
unequal to zero.
The length and weight of $\lambda$, denoted by 
$\ell(\lambda)$ and $\abs{\lambda}$, are the number and sum
of the non-zero $\lambda_i$ (called parts), respectively.
The unique partition of weight zero is denoted by $0$, and
the multiplicity of the part $i$ in the partition $\la$ is denoted
by $m_i(\lambda)$.

We identify a partition with its diagram or Ferrers graph in the
usual way, and, for example, the diagram of $\la=(6,3,3,1)$ is given by
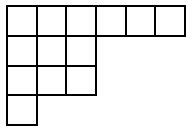
\begin{center}\yng(6,3,3,1)\end{center}

The conjugate $\lambda'$ of $\lambda$ is the partition obtained by 
reflecting the diagram of $\lambda$ in the main diagonal.
Hence $m_i(\lambda)=\la_i'-\la_{i+1}'$.

A standard statistic on partitions needed repeatedly is
\begin{equation*}
n(\la)=\sum_{i\geq 1} (i-1)\la_i=\sum_{i\geq 1}\binom{\la_i'}{2}.
\end{equation*}
We also need the usual scalar product $\pp{\la}{\mu}=\sum_{i\geq 1}\la_i\mu_i$
(which in the notation of \cite{Macdonald95} would be $\abs{\la\mu}$).
We will occasionally use this for more general sequences of integers, not
necessarily partitions.

If $\la$ and $\mu$ are two partions then $\mu\subset\la$ iff
$\la_i\geq \mu_i$ for all $i\geq 1$, i.e., the diagram of $\la$ contains
the diagram of $\mu$. If $\mu\subset\la$ then the skew-diagram $\la-\mu$
denotes the set-theoretic difference between $\la$ and $\mu$,
and $\abs{\la-\mu}=\abs{\la}-\abs{\mu}$.
For example, if $\la=(6,3,3,1)$ and $\mu=(4,3,1)$ then
the skew diagram $\la-\mu$ is given by the marked squares in
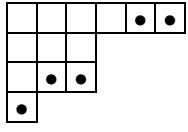
\begin{center}
\young(\hfil\hfil\hfil\hfil\vul\vul,\hfil\hfil\hfil,\hfil\vul\vul,\vul)
\end{center}
and $\abs{\la-\mu}=5$.

For $\theta=\la-\mu$ a skew diagram, its conjugate $\theta'=\la'-\mu'$ 
is the (skew) diagram obtained by reflecting $\theta$ in the main diagonal. 
Following \cite{Macdonald95} we define the components of $\theta$ and $\theta'$
by $\theta_i=\la_i-\mu_i$ and $\theta_i'=\la_i'-\mu_i'$.
Quite often we only require knowledge of
the sequence of components of a skew diagram $\theta$,
and by abuse of notation we will occasionally write 
$\theta=(\theta_1,\theta_2,\dots)$,
even though the components $\theta_i$ alone do not fix $\theta$.

A skew diagram $\theta$ is a horizontal strip if $\theta_i'\in\{0,1\}$,
i.e., if at most one square occurs in each column of $\theta$.
The skew diagram in the above example is a horizontal strip since
$\theta'=(1,1,1,0,1,1,0,0,\dots)$.

\medskip

Let $S_n$ be the symmetric group, $\Lambda_n=\Z[x_1,\dots,x_n]^{S_n}$ 
be the ring of symmetric polynomials in $n$ 
independent variables and $\Lambda$ the ring of symmetric functions
in countably many independent variables.

For $x=(x_1,\dots,x_n)$ and $\lambda$ a partition such that $\ell(\la)\leq n$
the Hall--Littlewood polynomials $P_{\lambda}(x;q)$ are defined by
\begin{equation}\label{Pdef} 
P_{\lambda}(x;q)
=\sum_{w\in S_n/S_n^{\la}} 
w\Bigl(x^{\la} \prod_{\la_i>\la_j}
\frac{x_i-qx_j}{x_i-x_j} \Bigr).
\end{equation}
Here $S_n^{\la}$ is the subgroup of $S_n$ consisting of the permutations 
that leave $\la$ invariant, and $w(f(x))=f(w(x))$.
When $\ell(\la)>n$, 
\begin{equation}\label{Pnul}
P_{\lambda}(x;q)=0.
\end{equation}

The Hall--Littlewood polynomials are symmetric polynomials in $x$,
homogeneous of degree $\abs{\la}$, with coefficients
in $\Z[q]$, and form a $\Z[q]$ basis of $\Lambda_n[q]$.
Thanks to the stability property $P_{\la}(x_1,\dots,x_n,0;q)=
P_{\la}(x_1,\dots,x_n;q)$ the Hall--Littlewood polynomials may be 
extended to the Hall--Littlewood functions in an infinite number of
variables $x_1,x_2,\dots$ in the usual way, to form a $\Z[q]$ basis of
$\Lambda[q]$.
The indeterminate $q$ in the Hall--Littlewood symmetric functions
serves as a parameter interpolating between the
Schur functions and monomial symmetric functions;
$P_{\la}(x;0)=s_{\la}(x)$ and $P_{\la}(x;1)=m_{\la}(x)$.

We will also need the symmetric functions $Q_{\lambda}(x;q)$ 
(also referred to as Hall-Littlewood functions) defined by
\begin{equation}\label{Qdef}
Q_{\la}(x;q)=b_{\la}(q)P_{\la}(x;q),
\end{equation}
where
\begin{equation*}
b_{\la}(q)=\prod_{i=1}^{\la_1}(q;q)_{m_i(\la)}.
\end{equation*}

We already mentioned the homogeneity of the
Hall--Littlewood functions;
\begin{equation}\label{hom}
P_{\la}(ax;q)=a^{\abs{\la}}P_{\la}(x;q),
\end{equation}
where $ax=(ax_1,ax_2,\dots)$.
Another useful result is the specialization
\begin{equation}\label{Pspec}
P_{\la}(1,q,\dots,q^{n-1};q)=
\frac{q^{n(\la)}(q;q)_n}
{(q;q)_{n-\ell(\la)}b_{\lambda}(q)},
\end{equation}
where $1/(q;q)_{-m}=0$ for $m$ a positive integer, so that
$P_{\la}(1,q,\dots,q^{n-1};q)=0$ if $\ell(\la)>n$ in accordance with
\eqref{Pnul}.
By \eqref{Qdef} this also implies the particularly simple
\begin{equation}\label{Qspec}
Q_{\la}(1,q,q^2,\dots;q)=q^{n(\la)}.
\end{equation}

The Cauchy identity for Hall--Littlewood functions states that
\begin{equation}\label{Cauchy}
\sum_{\la}P_{\la}(x;q)Q_{\la}(y;q)=\prod_{i,j\geq 1}\frac{1-qx_iy_j}{1-x_iy_j}.
\end{equation}
Taking $y_j=q^{j-1}$ for all $j\geq 1$ and using the specialization
\eqref{Qspec} yields
\begin{equation}\label{nPsum}
\sum_{\la}q^{n(\la)}P_{\la}(x;q)=\prod_{i\geq 1}\frac{1}{1-x_i}.
\end{equation}
We remark that this is the A$_1$ analogue of Theorem~\ref{mainthm},
providing an evaluation for the sum \eqref{Ansum} when $n=1$.

The skew Hall--Littlewood functions $P_{\la/\mu}$ and $Q_{\la/\mu}$
are defined by
\begin{equation}\label{Pskew}
P_{\la}(x,y;q)=\sum_{\mu}P_{\la/\mu}(x;q)P_{\mu}(y;q)
\end{equation}
and 
\begin{equation*}
Q_{\la}(x,y;q)=\sum_{\mu}Q_{\la/\mu}(x;q)Q_{\mu}(y;q),
\end{equation*}
so that 
\begin{equation}\label{QPskew}
Q_{\la/\mu}(x;q)=\frac{b_{\la}(q)}{b_{\mu}(q)}P_{\la/\mu}(x;q).
\end{equation}
An important property is that $P_{\la/\mu}$ is zero if $\mu\not\subset\la$.
Some trivial instances of the skew functions are given by
$P_{\la/0}=P_{\la}$ and $P_{\la/\la}=1$. 
By \eqref{QPskew} similar statements apply to $Q_{\la/\mu}$.

The Cauchy identity \eqref{Cauchy} can be generalized to
the skew case as \cite[Lemma 3.1]{Stembridge90}
\begin{equation}\label{SCauchy}
\sum_{\la}P_{\la/\mu}(x;q)Q_{\la/\nu}(y;q)=
\sum_{\la}P_{\nu/\la}(x;q)Q_{\mu/\la}(y;q)
\prod_{i,j\geq 1}\frac{1-qx_iy_j}{1-x_iy_j}.
\end{equation}
Taking $\nu=0$ and specializing $y_j=q^{j-1}$ for all $j\geq 1$ 
extends \eqref{nPsum} to 
\begin{equation}\label{SQspec}
\sum_{\la}q^{n(\la)}P_{\la/\mu}(x;q)=q^{n(\mu)}
\prod_{i\geq 1}\frac{1}{1-x_i}.
\end{equation}

We conclude our introduction of the Hall--Littlewood functions
with the following two important definitions.
Let $\la\supset\mu$ be partitions such that $\theta=\la-\mu$
is a horizontal strip, i.e., $\theta'_i\in\{0,1\}$. 
Let $I$ be the set of integers $i\geq 1$ such that $\theta'_i=1$ and
$\theta'_{i+1}=0$. Then
\begin{equation*}
\phi_{\la/\mu}(q)=\prod_{i\in I}(1-q^{m_i(\la)}).
\end{equation*}
Similarly, let $J$ be the set of integers $j\geq 1$ such that 
$\theta'_j=0$ and $\theta'_{j+1}=1$. Then
\begin{equation*}
\psi_{\la/\mu}(q)=\prod_{j\in J}(1-q^{m_j(\mu)}).
\end{equation*}
For example, if $\la=(5,3,2,2)$ and $\mu=(3,3,2)$ then $\theta$ is a horizontal
strip and $\theta'=(1,1,0,1,1,0,0,\dots)$. Hence $I=\{2,5\}$ and $J=\{3\}$,
leading to
\begin{equation*}
\phi_{\la/\mu}(q)=(1-q^{m_2(\la)})(1-q^{m_5(\la)})=(1-q^2)(1-q)
\end{equation*}
and
\begin{equation*}
\psi_{\la/\mu}(q)=(1-q^{m_3(\mu)})=(1-q^2).
\end{equation*}

The skew Hall--Littlewood functions $Q_{\la/\mu}(x;q)$ and $P_{\la/\mu}(x;q)$ 
can be expressed in terms of $\phi_{\la/\mu}(q)$ and $\psi_{\la/\mu}(q)$
\cite[p. 229]{Macdonald95}. For our purposes we only require
a special instance of this result corresponding to the case that 
$x$ represents a single variable. Then
\begin{subequations}\label{defphipsi}
\begin{equation}
Q_{\la/\mu}(x;q)
=\begin{cases}
\phi_{\la/\mu}(q)x^{\abs{\la-\mu}} &
\text{if $\la-\mu$ is a horizontal strip,} \\
0 &\text{otherwise}
\end{cases}
\end{equation}
and
\begin{equation}\label{P1var}
P_{\la/\mu}(x;q)
=\begin{cases}
\psi_{\la/\mu}(q)x^{\abs{\la-\mu}} &
\text{if $\la-\mu$ is a horizontal strip,} \\
0 &\text{otherwise.}
\end{cases}
\end{equation}
\end{subequations}

\section{Consequences of Theorem~\ref{mainthm}}\label{secC}
Before we present a proof of Theorem~\ref{mainthm}
we will establish some simple corollaries of the
A$_2$ sum for Hall--Littlewood functions.

We begin by noting that \eqref{Pid} simplifies to \eqref{nPsum}
when all components of $y$ are set to zero.
Our first corollary of Theorem~\ref{mainthm} corresponds to
a slight generalization that also includes \eqref{SQspec}.
\begin{corollary}\label{Cor}
For $\nu$ a partition,
\begin{multline}\label{Pidskew}
\sum_{\la,\mu} q^{n(\la)+n(\mu)-\pp{\la'}{\mu'}} P_{\la/\nu}(x;q)P_{\mu}(y;q) \\
=\sum_{\la} q^{n(\la)+n(\nu)-\pp{\la'}{\nu'}} P_{\la}(y;q)
\prod_{i\geq 1}\frac{1}{1-x_i}
\prod_{i,j\geq 1}\frac{1-x_iy_j}{1-q^{-1}x_iy_j}.
\end{multline}
\end{corollary}
When $\nu=0$ the sum over $\la$ on the right may be performed
by \eqref{nPsum} and one recovers \eqref{Pid}.

\begin{proof}[Proof of Corollary \ref{Cor}]
Multiplying both sides of \eqref{Pidskew} by
$P_{\nu}(z;q)$ and summing over $\nu$ gives
\begin{multline*}
\sum_{\la,\mu} q^{n(\la)+n(\mu)-\pp{\la'}{\mu'}} P_{\la}(x,z;q)P_{\mu}(y;q) \\
=\sum_{\la,\nu} q^{n(\la)+n(\nu)-\pp{\la'}{\nu'}} P_{\la}(y;q)P_{\nu}(z;q)
\prod_{i\geq 1}\frac{1}{1-x_i}
\prod_{i,j\geq 1}\frac{1-x_iy_j}{1-q^{-1}x_iy_j},
\end{multline*}
where on the left we have used \eqref{Pskew}.
The truth of this identity is readily verified upon noting that
both sides can be summed by \eqref{Pid}.
Since the $P_{\nu}(z;q)$ form a basis of $\Lambda[q]$ 
the identity \eqref{Pidskew} itself must be true.
\end{proof}

It is suggestive that a yet more general symmetric expansion should
hold for
\begin{equation*}
\sum_{\la,\mu} q^{n(\la)+n(\mu)-\pp{\la'}{\mu'}} 
P_{\la/\nu}(x;q)P_{\mu/\eta}(y;q),
\end{equation*}
but we were only able to obtain the following asymmetric sum.
\begin{corollary}\label{Cor2}
For $\nu$ and $\eta$ partitions,
\begin{multline*}
\sum_{\la,\mu} q^{n(\la)+n(\mu)-\pp{\la'}{\mu'}} 
P_{\la/\nu}(x;q)P_{\mu/\eta}(y;q) \\
=\sum_{\la,\mu} q^{n(\la)+n(\nu)-\pp{\la'}{\nu'}} 
Q_{\eta/\mu}(x/q;q)P_{\la/\mu}(y;q)
\prod_{i\geq 1}\frac{1}{1-x_i}
\prod_{i,j\geq 1}\frac{1-x_iy_j}{1-q^{-1}x_iy_j}.
\end{multline*}
\end{corollary}
When all the $y_i$ are set to zero this yields (after a change
of variables)
\begin{multline}\label{genx}
\sum_{\nu} q^{n(\la)+n(\nu)-\pp{\la'}{\nu'}} P_{\nu/\mu}(x;q) \\
=\sum_{\nu} q^{n(\mu)+n(\nu)-\pp{\mu'}{\nu'}} Q_{\la/\nu}(x/q;q)
\prod_{i\geq 1} \frac{1}{1-x_i}.
\end{multline}
The case when $x$ represents a single variable will
play an important role in the proof of Theorem~\ref{mainthm}.

\begin{proof}[Proof of Corollary~\ref{Cor2}]
After multiplying both sides by $P_{\nu}(z;q)P_{\eta}(w;q)$ 
and summing over $\nu$ and $\eta$ we get
\begin{multline*}
\sum_{\la,\mu} q^{n(\la)+n(\mu)-\pp{\la'}{\mu'}} 
P_{\la}(x,z;q)P_{\mu}(y,w;q) \\
=\sum_{\la,\mu,\nu,\eta} q^{n(\la)+n(\nu)-\pp{\la'}{\nu'}} 
Q_{\eta/\mu}(x/q;q)P_{\la/\mu}(y;q)P_{\nu}(z;q)P_{\eta}(w;q)
\\ \times
\prod_{i\geq 1}\frac{1}{1-x_i}
\prod_{i,j\geq 1}\frac{1-x_iy_j}{1-q^{-1}x_iy_j}.
\end{multline*}
By the Cauchy identity for skew Hall--Littlewood functions 
\eqref{SCauchy} the sum over $\eta$ on the right simplifies to
\begin{equation*}
P_{\mu}(w;q)\prod_{i,j\geq 1}\frac{1-x_iw_j}{1-q^{-1}x_iw_j}.
\end{equation*}
This then allows for the sum over $\mu$ to be carried out 
using \eqref{Pskew}, yielding
\begin{multline*}
\sum_{\la,\mu} q^{n(\la)+n(\mu)-\pp{\la'}{\mu'}} 
P_{\la}(x,z;q)P_{\mu}(y,w;q) \\
=\sum_{\la,\nu} q^{n(\la)+n(\nu)-\pp{\la'}{\nu'}} 
P_{\la}(y,w;q)P_{\nu}(z;q) \\ \times
\prod_{i\geq 1}\frac{1}{1-x_i}
\prod_{i,j\geq 1}\frac{(1-x_iy_j)(1-x_iw_j)}
{(1-q^{-1}x_iy_j)(1-q^{-1}x_iw_j)}.
\end{multline*}
The rest again follows from \eqref{Pid}.
\end{proof}

As a third corollary we can include a linear term in
the exponent of $q$ in \eqref{nPsum}.
\begin{corollary}\label{Cor3}
For $j$ a non-negative integer,
\begin{equation*}
\sum_{\la} q^{n(\la)-\sum_{l=1}^j \la'_l} 
P_{\la}(x;q)
=\Bigl(1+(1-q)\sum_{k=1}^j q^{-k} P_{(k)}(x;q)\Bigr)
\prod_{i\geq 1}\frac{1}{1-x_i}.
\end{equation*}
\end{corollary}

This result implies some nice $q$-series identities.
By specializing $x_i=zq^i$ for $1\leq i\leq n$ and $x_i=0$ for $i>n$,
and using \eqref{hom} and \eqref{Pspec} we find
\begin{align*}
\sum_{\la} \frac{z^{\abs{\la}}q^{\pp{\la'}{\la'}-\sum_{l=1}^j \la'_l}(q;q)_n}
{(q;q)_{n-\ell(\la)}b_{\la}(q)}&=
\frac{1}{(zq;q)_n}
\Bigl(1+(1-q^n)\sum_{k=1}^j z^k\Bigr) \\
&=\frac{1}{(z;q)_n}-(1-q^n)\frac{z^{j+1}}{(z;q)_{n+1}},
\end{align*}
where we have used that $2n(\la)+\abs{\la}=\pp{\la'}{\la'}$.
By the $q$-binomial theorem
\begin{equation*}
\sum_{j=0}^{\infty} z^j \qbin{n+j-1}{j}=\frac{1}{(z;q)_n}
\end{equation*}
the coefficient of $z^k$ can easily be found as
\begin{equation*}
\sum_{\lambda\vdash k}
\frac{q^{\pp{\la}{\la}-\sum_{l=1}^j \la_l}(q;q)_n}
{(q;q)_{n-\la_1(\la)}b_{\la'}(q)}
=\qbin{n+k-1}{k}-(1-q^n)\qbin{n+k-j-1}{k-j-1}
\end{equation*}
for $0\leq j\leq k$. Here we have changed the summation index $\la$ by
its conjugate. For $j=k$ or (after simplifying the right) $j=0$
this is a well-known $q$-series identity of
Hall \cite{Hall38}, see also \cite{Macdonald89,SZ89,Stembridge90}.
Letting $n$ tend to infinity finally gives
\begin{equation*}
\sum_{\lambda\vdash k}
\frac{q^{\pp{\la}{\la}-\sum_{l=1}^j \la_j}}{b_{\la'}(q)}
=\frac{1}{(q;q)_k}-\frac{1}{(q;q)_{k-j-1}}.
\end{equation*}

\begin{proof}[Proof of Corollary~\ref{Cor3}]
Equation \eqref{Pid} with $y_1=z$ and $y_i=0$ for $i\geq 2$ yields
\begin{equation*}
\sum_{j=0}^{\infty} z^j \sum_{\la}
q^{n(\la)-\sum_{l=1}^j \la'_l} 
P_{\la}(x;q) \\
=\frac{1}{1-z}\prod_{i\geq 1}\frac{1}{1-x_i}
\prod_{i\geq 1}\frac{1-x_iz}{1-q^{-1}x_iz}.
\end{equation*}
From the Cauchy identity \eqref{Cauchy} it follows that
the last product on the right can be expanded as
\begin{equation*}
\sum_{k=0}^{\infty} Q_{(k)}(z/q;q)P_{(k)}(x;q)=
1+(1-q)\sum_{k=0}^{\infty} P_{(k)}(x;q)(z/q)^k.
\end{equation*}
Then equating coefficients of $z^j$ leads to the desired result.
\end{proof}

Finally we come to what is by far the most important corollary of
Theorem~\ref{mainthm}.
Let $(a;q)_0=1$, $(a;q)_n=\prod_{i=1}^n (1-aq^{i-1})$ and
$(a_1,\dots,a_k;q)_n=(a_1;q)_n\cdots(a_k;q)_n$.
\begin{corollary}\label{Cor4}
There holds
\begin{equation}\label{PPs}
\sum_{\la,\mu} \frac{a^{\abs{\la}}b^{\abs{\mu}}
q^{\pp{\la'}{\la'}+\pp{\mu'}{\mu'}-\pp{\la'}{\mu'}}}
{(q;q)_{n-\ell(\la)}(q;q)_{m-\ell(\mu)}b_{\la}(q)b_{\mu}(q)} \\
=\frac{(abq;q)_{n+m}}{(q,aq,abq;q)_n(q,bq,abq;q)_m}
\end{equation}
\end{corollary}

\begin{proof}
Taking $x_i=aq^i$ for $1\leq i\leq n$, $x_i=0$ for $i>n$,
$y_j=bq^j$ for $1\leq j\leq m$ and $y_j=0$ for $j>m$, 
using the homogeneity \eqref{hom} and specialization \eqref{Pspec},
and noting that $2n(\la)+\abs{\la}=\pp{\la'}{\la'}$,
we obtain \eqref{PPs}.
\end{proof}

In Section~\ref{secRR} we will show how Corollary~\ref{Cor4} relates
to the A$_2$ Rogers--Ramanujan identity \eqref{RRA2}.
For now let us remark that \eqref{PPs} is a bounded
version of the A$_2$ case of the following identity for
the A$_n$ root system due to Hua \cite{Hua00} (and
corrected in \cite{Fulman01}):
\begin{equation}\label{Hua}
\sum_{\la^{(1)},\dots,\la^{(n)}}
\frac{q^{\frac{1}{2}\sum_{i,j=1}^n C_{ij}\pp{{\la^{(i)}}'}{{\la^{(j)}}'}}
\prod_{i=1}^n a_i^{\abs{\la^{(i)}}}}
{\prod_{i=1}^n b_{\la^{(i)}}(q)}
=\prod_{\alpha\in\Delta_{+}}\frac{1}{(a^{\alpha}q;q)_{\infty}}.
\end{equation}
Here $C_{ij}=2\delta_{i,i}-\delta_{i,j-1}-\delta_{i,j+1}$ is 
the $(i,j)$ entry of the A$_n$ Cartan matrix and
$\Delta_{+}$ is the set of positive roots of A$_n$, i.e.,
the set (of cardinality $\binom{n+1}{2}$) of
roots of the form $\alpha_i+\alpha_{i+1}+\cdots+\alpha_j$ 
with $1\leq i\leq j\leq n$, where
$\alpha_1,\dots,\alpha_n$ are the simple roots of A$_n$.
Furthermore, if $\alpha=\alpha_i+\alpha_{i+1}+\cdots+\alpha_j$ then
$a^{\alpha}=a_ia_{i+1}\cdots a_j$.

\section{Proof of Theorem~\ref{mainthm}}\label{secP}
Throughout this section $z$ represents a single variable.

To establish \eqref{Pid} it is enough to show its truth
for $x=(x_1,\dots,x_n)$ and $y=(y_1,\dots,y_m)$, 
and by induction on $m$ it then easily follows that we
only need to prove
\begin{multline}\label{rec}
\sum_{\la,\mu} q^{n(\la)+n(\mu)-\pp{\la'}{\mu'}} P_{\la}(x;q)P_{\mu}(y,z;q) \\
=\frac{1}{1-z}
\prod_{i=1}^n\frac{1-zx_i}{1-q^{-1}zx_i}
\sum_{\la,\mu} q^{n(\la)+n(\mu)-\pp{\la'}{\mu'}} P_{\la}(x;q)P_{\mu}(y;q),
\end{multline}
where we have replaced $y_{m+1}$ by $z$.

If on the left we replace $\mu$ by $\nu$ and use \eqref{Pskew}
(with $\la\to\nu$ and $x\to z$) we get
\begin{equation*}
\text{LHS}\eqref{rec}=
\sum_{\la,\mu,\nu} 
q^{n(\la)+n(\nu)-\pp{\la'}{\nu'}} P_{\la}(x;q) P_{\mu}(y;q) P_{\nu/\mu}(z;q).
\end{equation*}
From \eqref{SCauchy} with $\mu=0$, $x=(x_1,\dots,x_n)$ and $y\to z/q$
it follows that
\begin{equation*}
P_{\nu}(x;q) \prod_{i=1}^n\frac{1-zx_i}{1-q^{-1}zx_i}
=\sum_{\la}Q_{\la/\nu}(z/q;q)P_{\la}(x;q).
\end{equation*}
Using this on the right of \eqref{rec} with $\la$ replaced by $\nu$ 
yields
\begin{equation*}
\text{RHS}\eqref{rec}=\frac{1}{1-z}
\sum_{\la,\mu,\nu}
q^{n(\mu)+n(\nu)-\pp{\mu'}{\nu'}} P_{\la}(x;q)P_{\mu}(y;q)Q_{\la/\nu}(z/q;q).
\end{equation*}
Therefore, by equating coefficients of $P_{\la}(x;q)P_{\mu}(y;q)$ 
we find that the problem of proving \eqref{Pid} boils down to showing that
\begin{equation*}
\sum_{\nu}
q^{n(\la)+n(\nu)-\pp{\la'}{\nu'}} P_{\nu/\mu}(z;q)
=\frac{1}{1-z}
\sum_{\nu}
q^{n(\mu)+n(\nu)-\pp{\mu'}{\nu'}} Q_{\la/\nu}(z/q;q),
\end{equation*}
which is \eqref{genx} with $x\to z$.

Next we use \eqref{defphipsi} to arrive at the equivalent but
more combinatorial statement that
\begin{multline}\label{psiphi}
\sum_{\substack{\nu\supset\mu \\ \nu-\mu \textup{ hor. strip}}}
q^{n(\la)+n(\nu)-\pp{\la'}{\nu'}} z^{\abs{\nu-\mu}}\psi_{\nu/\mu}(q) \\
=\frac{1}{1-z}
\sum_{\substack{\nu\subset\la \\ \la-\nu \textup{ hor. strip}}}
q^{n(\mu)+n(\nu)-\pp{\mu'}{\nu'}} (z/q)^{\abs{\la-\nu}}\phi_{\la/\nu}(q).
\end{multline}
This identity is reminiscient of the well-known but much simpler
\cite[Eq. (1); p. 230]{Macdonald95}
\begin{equation*}
\sum_{\substack{\mu\supset\la \\ \mu-\la \textup{ hor. strip}}}
z^{\abs{\mu-\la}}\psi_{\mu/\la}(q)
=\frac{1}{1-z}
\sum_{\substack{\mu\subset\la \\ \la-\mu \textup{ hor. strip}}}
z^{\abs{\la-\mu}}\phi_{\la/\mu}(q).
\end{equation*}

To make further progress we need a lemma.
\begin{lemma}\label{lem}
For $k$ a positive integer let $\omega=(\omega_1,\dots,\omega_k)\in
\{0,1\}^k$, and let $J=J(\omega)$ be the set of 
integers $j$ such that $\omega_j=0$ and $\omega_{j+1}=1$. 
For $\lambda\supset\mu$ partitions let $\theta'=\lambda'-\mu'$
be a skew diagram. Then
\begin{multline*}
\sum_{\substack{\la\supset\mu \\ \la-\mu \textup{ hor. strip} 
\\ \theta'_i=\omega_i,~i\in\{1,\dots,k\}}}
q^{n(\la)} z^{\abs{\la-\mu}}\psi_{\la/\mu}(q) \\
=\frac{q^{n(\mu)+\pp{\mu'}{\omega}} z^{\abs{\omega}}}{1-z}
(1-z(1-\omega_k)q^{\mu'_k})
\prod_{j\in J}(1-q^{m_j(\mu)}).
\end{multline*}
\end{lemma}
The restriction $\theta'_i=\omega_i$ for $i\in\{1,\dots,k\}$
in the sum over $\la$ on the left means that the first $k$ parts of
$\la'$ are fixed. The remaining parts are
free subject only to the condition that $\la-\mu$
is a horizontal strip, i.e., that $\la_i'-\mu_i'\in\{0,1\}$.

\begin{proof}
From \eqref{SQspec} with $x\to z$ combined with \eqref{P1var} we have
\begin{equation}\label{knul}
\sum_{\substack{\la\supset\mu \\ \la-\mu \textup{ hor. strip}}}
q^{n(\la)} z^{\abs{\la-\mu}}\psi_{\la/\mu}(q)
=\frac{q^{n(\mu)}}{1-z}.
\end{equation}
We will use this to first prove the lemma for $k=1$.
When $k=1$ and $\omega_1=1$ we need to show that
\begin{equation}\label{komegaeen}
\sum_{\substack{\la\supset\mu \\ \la-\mu \textup{ hor. strip} 
\\ \theta'_1=1}}
q^{n(\la)} z^{\abs{\la-\mu}}\psi_{\la/\mu}(q)
=\frac{q^{n(\mu)+\mu'_1} z}{1-z}.
\end{equation}
Now let $\lab$ and $\mub$ be the partitions obtained
from $\la$ and $\mu$ by removal of the first column of
their respective diagrams;
$\lab=(\la'_2,\la'_3,\dots)'$ and $\mub=(\mu'_2,\mu'_3,\dots)'$.
Since $\theta'_1=1$ we have 
$\psi_{\la/\mu}(q)=\psi_{\lab/\mub}(q)$, 
$\abs{\la-\mu}=\abs{\lab-\mub}+1$ and $\la'_1=\mu'_1+1$ so that
$n(\la)=n(\lab)+\binom{\mu'_1+1}{2}$.
Hence
\begin{align*}
\sum_{\substack{\la\supset\mu \\ \la-\mu \textup{ hor. strip} 
\\ \theta'_1=1}} q^{n(\la)} z^{\abs{\la-\mu}}\psi_{\la/\mu}(q) 
&=z q^{\binom{\mu'_1+1}{2}}
\sum_{\substack{\lab\supset\mub \\ \lab-\mub \textup{ hor. strip}}}
q^{n(\lab)} z^{\abs{\lab-\mub}}\psi_{\lab/\mub}(q) \\
&=z q^{\binom{\mu'_1+1}{2}}
\times \frac{q^{n(\mub)}}{1-z} \\
&=\frac{q^{n(\mu)+\mu'_1} z}{1-z},
\end{align*}
where the second equality follows from \eqref{knul}.

When $k=1$ and $\omega_1=0$ we need to show that
\begin{equation*}
\sum_{\substack{\la\supset\mu \\ \la-\mu \textup{ hor. strip} 
\\ \theta'_1=0}}
q^{n(\la)} z^{\abs{\la-\mu}}\psi_{\la/\mu}(q)
=\frac{q^{n(\mu)}}{1-z}(1-zq^{\mu'_1}).
\end{equation*}
This time we cannot simply relate $\psi_{\la/\mu}(q)$ to $\psi_{\lab/\mub}(q)$, 
but by inclusion-exclusion we have
\begin{align*}
\sum_{\substack{\la\supset\mu \\ \la-\mu \textup{ hor. strip} 
\\ \theta'_1=0}}&
q^{n(\la)} z^{\abs{\la-\mu}}\psi_{\la/\mu}(q) \\
&=\Biggl(\quad\sum_{\substack{\la\supset\mu \\ \la-\mu \textup{ hor. strip}}}
-\sum_{\substack{\la\supset\mu \\ \la-\mu \textup{ hor. strip} 
\\ \theta'_1=1}}\quad\Biggr)
q^{n(\la)} z^{\abs{\la-\mu}}\psi_{\la/\mu}(q) \\
&=\frac{q^{n(\mu)}}{1-z}-\frac{q^{n(\mu)+\mu'_1} z}{1-z} \\
&=\frac{q^{n(\mu)}}{1-z}(1-zq^{\mu'_1}),
\end{align*}
where the second equality follows from \eqref{knul} and \eqref{komegaeen}.

The remainder of the proof proceeds by induction on $k$.
Let us assume the lemma to be true for all $1\leq k\leq K-1$ with $K\geq 2$,
and use this to show its truth for $k=K$.
To do so we need to again distinguish two cases:
$\omega=(\omega_1,\dots,\omega_K)$ with (i)
$\omega_1=1$ or $\omega_1=\omega_2=0$ and (ii) $\omega_1=0$ and $\omega_2=1$.

First consider (i) and attach the same meaning to $\lab$ and $\mub$ as before.
We also set $\omegab=(\omegab_1,\dots,\omegab_{K-1})=(\omega_2,\dots,\omega_K)$
and $\thetab_i=\theta_{i+1}$. Then, since $\omega_1=1$ or $\omega_1=\omega_2=0$,
\begin{equation*}
\psi_{\la/\mu}(q)=\psi_{\lab/\mub}(q) \quad \text{and}\quad
\prod_{j\in J(\omega)}(1-q^{m_j(\mu)})
=\prod_{j\in J(\omegab)}(1-q^{m_j(\mub)}).
\end{equation*}
Moreover, 
$\abs{\la-\mu}=\abs{\lab-\mub}+\omega_1$ and $\la'_1=\mu'_1+\omega_1$ so that
$n(\la)=n(\lab)+\binom{\mu'_1+\omega_1}{2}$.
Therefore
\begin{align*}
&\sum_{\substack{\la\supset\mu \\ \la-\mu \textup{ hor. strip} 
\\ \theta'_i=\omega_i,~i\in\{1,\dots,K\}}}
q^{n(\la)} z^{\abs{\la-\mu}}\psi_{\la/\mu}(q) \\
&\quad=z^{\omega_1}q^{\binom{\mu'_1+\omega_1}{2}}
\sum_{\substack{\lab\supset\mub \\ \lab-\mub \textup{ hor. strip} 
\\ \thetab'_i=\omegab_i,~i\in\{1,\dots,K-1\}}}
q^{n(\lab)} z^{\abs{\lab-\mub}}\psi_{\lab/\mub}(q) \\[2mm]
&\quad=z^{\omega_1}q^{\binom{\mu'_1+\omega_1}{2}}\times
\frac{q^{n(\mub)+\pp{\mub'}{\omegab}} z^{\abs{\omegab}}}{1-z}
(1-z(1-\omegab_{K-1})q^{\mub'_{K-1}})
\prod_{j\in J(\omegab)}(1-q^{m_j(\mub)}) \\
&\quad=\frac{q^{n(\mu)+\pp{\mu'}{\omega}} z^{\abs{\omega}}}{1-z}
(1-z(1-\omega_K)q^{\mu'_K})
\prod_{j\in J(\omega)}(1-q^{m_j(\mu)}).
\end{align*}

In the case of (ii) the proof requires only minor changes, and this time
we need
\begin{gather*}
\psi_{\la/\mu}(q)=\psi_{\lab/\mub}(q)(1-q^{m_1(\mu)}), \\
\intertext{and}
\prod_{j\in J(\omega)}(1-q^{m_j(\mu)})
=(1-q^{m_1(\mu)})\prod_{j\in J(\omegab)}(1-q^{m_j(\mub)}).
\end{gather*}
(Note that both sides of the first of these equations vanish 
if $m_1(\mu)=0$ as it should.
Indeed, if $\mu'_1=\mu'_2$ there is no partition $\la\supset\mu$ such that
$\theta'=\la'-\mu'=(0,1,\omega_3,\dots,\omega_K)$ since it would require
that $\la'_1<\la_2$.)
\end{proof}

In view of Lemma~\ref{lem} it is natural to rewrite the
left side of \eqref{psiphi} as
\begin{equation*}
\text{LHS}\eqref{psiphi}=
\sum_{\omega\in\{0,1\}^{\la_1}}
\sum_{\substack{\nu\supset\mu \\ \nu-\mu \text{ hor. strip} 
\\ \theta'_i=\omega_i,~i\in\{1,\dots,\la_1\}}}
q^{n(\la)+n(\nu)-\pp{\la'}{\mu'}-\pp{\la'}{\omega}} 
z^{\abs{\nu-\mu}}\psi_{\nu/\mu}(q),
\end{equation*}
where $\theta=\nu-\mu$, and where we have used that 
$\theta'_i\in\{0,1\}$ as follows from the fact that
$\nu-\mu$ is a horizontal strip.

Now the sum over $\nu$ can be performed by application
of Lemma~\ref{lem} with $\la\to\nu$ and $k\to \la_1$, resulting in
\begin{multline*}
\text{LHS}\eqref{psiphi}=
\frac{q^{n(\la)+n(\mu)-\pp{\la'}{\mu'}}}{1-z}
\sum_{\omega\in\{0,1\}^{\la_1}}
q^{\pp{\mu'}{\omega}-\pp{\la'}{\omega}} 
z^{\abs{\omega}} \\ \times
(1-z(1-\omega_{\la_1})q^{\mu'_{\la_1}})
\prod_{j\in J}(1-q^{m_j(\mu)})
\end{multline*}
with $J=J(\omega)\subset \{1,\dots,\la_1-1\}$ the set of
integers $j$ such that $\omega_j<\omega_{j+1}$.

For the right-hand side of \eqref{psiphi} we introduce the notation
$\tau_i=\la_i'-\nu_i'$, so that the sum over $\nu$ 
can be rewritten as a sum over $\tau\in\{0,1\}^{\la_1}$.
Using that 
\begin{equation*}
n(\nu)=\sum_{i=1}^{\la_1}\binom{\nu_i'}{2}=
\sum_{i=1}^{\la_1}\binom{\la_i'-\tau_i}{2}=
n(\la)-\pp{\la'}{\tau}+\abs{\tau}
\end{equation*}
this yields
\begin{equation*}
\text{RHS}\eqref{psiphi}=
\frac{q^{n(\la)+n(\mu)-\pp{\la'}{\mu'}}}{1-z}
\sum_{\tau\in\{0,1\}^{\la_1}}
q^{\pp{\mu'}{\tau}-\pp{\la'}{\tau}} 
z^{\abs{\tau}}
\prod_{i\in I}(1-q^{m_i(\la)}),
\end{equation*}
with $I=I(\tau)\subset \{1,\dots,\la_1\}$ the set of integers
$i$ such that $\tau_i>\tau_{i+1}$ (with the convention that 
$\la_1 \in I$ if $\tau_{\la_1}=1$).

Equating the above two results for the respective sides of \eqref{psiphi}
gives
\begin{multline*}
\sum_{\omega\in\{0,1\}^{\la_1}}
q^{\pp{\mu'}{\omega}-\pp{\la'}{\omega}} 
z^{\abs{\omega}} (1-z(1-\omega_{\la_1})q^{\mu'_{\la_1}})
\prod_{j\in J}(1-q^{m_j(\mu)}) \\
=\sum_{\tau\in\{0,1\}^{\la_1}}
q^{\pp{\mu'}{\tau}-\pp{\la'}{\tau}} z^{\abs{\tau}}
\prod_{i\in I}(1-q^{m_i(\la)}).
\end{multline*}
Using that $m_i(\la)=\la'_i-\la'_{i+1}$ it is not hard to see
that this is the 
\begin{equation*}
k\to \la_1,\quad b_{k+1}\to 1,\quad
a_i\to zq^{\mu_i'},\quad b_i\to q^{\la_i'},\quad i\in\{1,\dots,\la_1\}
\end{equation*}
specialization of the more general
\begin{multline*}
\sum_{\omega\in\{0,1\}^k} (a/b)^{\omega}
(1-(1-\omega_k)a_k/b_{k+1})
\prod_{j\in J}(1-a_j/a_{j+1}) \\
=\sum_{\tau\in\{0,1\}^k}
(a/b)^{\tau} \prod_{i\in I}(1-b_i/b_{i+1}),
\end{multline*}
where $(a/b)^{\omega}=\prod_{i=1}^k (a_i/b_i)^{\omega_i}$
and $(a/b)^{\tau}=\prod_{i=1}^k (a_i/b_i)^{\tau_i}$.
Obviously, the set $J\subset \{1,\dots,k-1\}$ should now be defined as
the set of integers $j$ such that $\omega_j<\omega_{j+1}$
and the the set $I\subset \{1,\dots,k\}$ as the set of integers
$i$ such that $\tau_i>\tau_{i+1}$ (with the convention that
$k\in I$ if $\tau_k=1$).

Next we split both sides into the sum of two terms as follows:
\begin{multline*}
\Biggl(\;\sum_{\omega\in\{0,1\}^k} 
-(a_k/b_{k+1})\sum_{\substack{\omega\in\{0,1\}^k \\ \omega_k=0}}\;
\Biggr) (a/b)^{\omega}
\prod_{j\in J}(1-a_j/a_{j+1}) \\
=\Biggr(\;\sum_{\tau\in\{0,1\}^k}
-(b_k/b_{k+1})\sum_{\substack{\tau\in\{0,1\}^k \\ \tau_k=1}}\;
\Biggl)
(a/b)^{\tau} \prod_{\substack{i\in I \\ i\neq k}}(1-b_i/b_{i+1}).
\end{multline*}
Equating the first sum on the left with the first sum on the right
yields
\begin{equation}\label{ab}
\sum_{\omega\in\{0,1\}^k} (a/b)^{\omega}
\prod_{j\in J}(1-a_j/a_{j+1}) \\
=\sum_{\tau\in\{0,1\}^k}
(a/b)^{\tau}\prod_{\substack{i\in I \\ i\neq k}}(1-b_i/b_{i+1}).
\end{equation}
If we equate the second sum on the left with the second sum on the right
and use that $k-1\not\in J(\omega)$ if $\omega_k=0$ and
$k-1\not\in I(\tau)$ if $\tau_k=1$,
we obtain $(a_k/b_{k+1})(\eqref{ab}_{k\to k-1})$.

Slightly changing our earlier convention we thus need to prove that
\begin{equation}\label{ab2}
\sum_{\omega\in\{0,1\}^k} (a/b)^{\omega}
\prod_{j\in J}(1-a_j/a_{j+1}) \\
=\sum_{\tau\in\{0,1\}^k}
(a/b)^{\tau} \prod_{i\in I}(1-b_i/b_{i+1}),
\end{equation}
where from now on $I\subset\{1,\dots,k-1\}$ denotes the set
of integers $i$ such that $\tau_i>\tau_{i+1}$
(so that no longer $k\in I$ if $\tau_k=1$).
It is not hard to see by multiplying out the 
respective products that boths sides yield
$((1+\sqrt{2})^{k+1}-(1-\sqrt{2})^{k+1})/(2\sqrt{2})$
terms. To see that the terms on the left and right
are in one-to-one correspondence we again resort to
induction. First, for $k=1$ it is readily checked that
both sides yield $1+a_1/b_1$. For $k=2$ we on the left get
\begin{equation*}
\underbrace{1}_{\omega=(0,0)}\;+\;
\underbrace{(a_1/b_1)}_{\omega=(1,0)}\;+\;
\underbrace{(a_2/b_2)(1-a_1/a_2)}_{\omega=(0,1)}\;+\;
\underbrace{(a_1a_2/b_1b_2)}_{\omega=(1,1)}
\end{equation*}
and on the right
\begin{equation*}
\underbrace{1}_{\tau=(0,0)}\;+\;
\underbrace{(a_1/b_1)(1-b_1/b_2)}_{\tau=(1,0)}\;+\;
\underbrace{(a_2/b_2)}_{\tau=(0,1)}\;+\;
\underbrace{(a_1a_2/b_1b_2)}_{\tau=(1,1)}
\end{equation*}
which both give
\begin{equation*}
1+a_1/b_1+a_2/b_2-a_1/b_2+a_1a_2/b_1b_2.
\end{equation*}
Let us now assume that \eqref{ab2} has been shown to be true for
$1\leq k\leq K-1$ with $K\geq 3$ and prove the case $k=K$.

On the left of \eqref{ab2} we split the sum over $\omega$ according to
\begin{equation*}
\sum_{\omega\in\{0,1\}^k}=
\sum_{\substack{\omega\in\{0,1\}^k \\ \omega_1=1}}+
\sum_{\substack{\omega\in\{0,1\}^k \\ \omega_1=\omega_2=0}}+
\sum_{\substack{\omega\in\{0,1\}^k \\ \omega_1=0,~\omega_2=1}}.
\end{equation*}
Defining $\bar{\omega}\in\{0,1\}^{k-1}$ and
$\bar{\bar{\omega}}\in\{0,1\}^{k-2}$ by 
$\bar{\omega}=(\omega_2,\dots,\omega_k)$ and
$\bar{\bar{\omega}}=(\omega_3,\dots,\omega_k)$, 
and also setting
and $\bar{a}_j=a_{j+1}$, $\bar{b}_j=b_{j+1}$,
and $\bar{\bar{a}}_j=a_{j+2}$, $\bar{\bar{b}}_j=b_{j+2}$,
this leads to
\begin{align*}
\text{LHS}\eqref{ab2}&=
(a_1/b_1)\sum_{\bar{\omega}\in\{0,1\}^{k-1}} (\bar{a}/\bar{b})^{\bar{\omega}}
\prod_{j\in J(\bar{\omega})}(1-\bar{a}_j/\bar{a}_{j+1}) \\
&\quad+\sum_{\substack{\bar{\omega}\in\{0,1\}^{k-1} \\ \bar{\omega}_1=0}} 
(\bar{a}/\bar{b})^{\bar{\omega}}
\prod_{j\in J(\bar{\omega})}(1-\bar{a}_j/\bar{a}_{j+1}) \\
&\quad+(1-a_1/a_2)\sum_{\substack{\bar{\omega}\in\{0,1\}^{k-1}\\ \bar{\omega}_1=1}} 
(\bar{a}/\bar{b})^{\bar{\omega}}
\prod_{j\in J(\bar{\omega})}(1-\bar{a}_j/\bar{a}_{j+1}) \\
&=(1+a_1/b_1)\sum_{\bar{\omega}\in\{0,1\}^{k-1}} (\bar{a}/\bar{b})^{\bar{\omega}}
\prod_{j\in J(\bar{\omega})}(1-\bar{a}_j/\bar{a}_{j+1}) \\
&\quad-(a_1/a_2)\sum_{\substack{\bar{\omega}\in\{0,1\}^{k-1}\\ \bar{\omega}_1=1}} 
(\bar{a}/\bar{b})^{\bar{\omega}}
\prod_{j\in J(\bar{\omega})}(1-\bar{a}_j/\bar{a}_{j+1}) \\
&=(1+a_1/b_1)\sum_{\bar{\omega}\in\{0,1\}^{k-1}} (\bar{a}/\bar{b})^{\bar{\omega}}
\prod_{j\in J(\bar{\omega})}(1-\bar{a}_j/\bar{a}_{j+1}) \\
&\quad-(a_1/b_2)\sum_{\bar{\bar{\omega}}\in\{0,1\}^{k-2}}
(\bar{\bar{a}}/\bar{\bar{b}})^{\bar{\bar{\omega}}}
\prod_{j\in J(\bar{\bar{\omega}})}(1-\bar{\bar{a}}_j/\bar{\bar{a}}_{j+1}).
\end{align*}

On the right of \eqref{ab2} we split the sum over $\tau$ according to
\begin{equation*}
\sum_{\tau\in\{0,1\}^k}=
\sum_{\substack{\tau\in\{0,1\}^k \\ \tau_1=0}}+
\sum_{\substack{\tau\in\{0,1\}^k \\ \tau_1=\tau_2=1}}+
\sum_{\substack{\tau\in\{0,1\}^k \\ \tau_1=1,~\tau_2=0}}.
\end{equation*}
Defining $\bar{\tau}\in\{0,1\}^{k-1}$ and
$\bar{\bar{\tau}}\in\{0,1\}^{k-2}$ by   
$\bar{\tau}=(\tau_2,\dots,\tau_k)$ and
$\bar{\bar{\tau}}=(\tau_3,\dots,\tau_k)$,
this yields
\begin{align*}
\text{RHS}\eqref{ab2}&=
\sum_{\bar{\tau}\in\{0,1\}^{k-1}} (\bar{a}/\bar{b})^{\bar{\tau}}
\prod_{j\in J(\bar{\tau})}(1-\bar{b}_j/\bar{b}_{j+1}) \\
&\quad+(a_1/b_1)\sum_{\substack{\bar{\tau}\in\{0,1\}^{k-1} \\ \bar{\tau}_1=1}} 
(\bar{a}/\bar{b})^{\bar{\tau}}
\prod_{j\in J(\bar{\tau})}(1-\bar{b}_j/\bar{b}_{j+1}) \\
&\quad+(a_1/b_1)(1-b_1/b_2)\sum_{\substack{\bar{\tau}\in\{0,1\}^{k-1}\\ 
\bar{\tau}_1=0}} (\bar{a}/\bar{b})^{\bar{\tau}}
\prod_{j\in J(\bar{\tau})}(1-\bar{b}_j/\bar{b}_{j+1}) \\
&=(1+a_1/b_1)\sum_{\bar{\tau}\in\{0,1\}^{k-1}} (\bar{a}/\bar{b})^{\bar{\tau}}
\prod_{j\in J(\bar{\tau})}(1-\bar{b}_j/\bar{b}_{j+1}) \\
&\quad-(a_1/b_2)\sum_{\substack{\bar{\tau}\in\{0,1\}^{k-1}\\ \bar{\tau}_1=0}} 
(\bar{a}/\bar{b})^{\bar{\tau}}
\prod_{j\in J(\bar{\tau})}(1-\bar{b}_j/\bar{b}_{j+1}) \\
&=(1+a_1/b_1)\sum_{\bar{\tau}\in\{0,1\}^{k-1}} (\bar{a}/\bar{b})^{\bar{\tau}}
\prod_{j\in J(\bar{\tau})}(1-\bar{b}_j/\bar{b}_{j+1}) \\
&\quad-(a_1/b_2)\sum_{\bar{\bar{\tau}}\in\{0,1\}^{k-2}}
(\bar{\bar{a}}/\bar{\bar{b}})^{\bar{\bar{\tau}}}
\prod_{j\in J(\bar{\bar{\tau}})}(1-\bar{\bar{b}}_j/\bar{\bar{b}}_{j+1}).
\end{align*}
By our induction hypothesis this equates with the 
previous expression for the left-hand side of 
\eqref{ab2}, completing the proof.

\section{Corollary~\ref{Cor4} and the A$_2$ Rogers--Ramanujan identities}
\label{secRR}

For $M=(M_1,\dots,M_n)$ with $M_i$ a non-negative integer,
and $C$ the A$_n$ Cartan matrix
we define the following bounded analogue of the sum in \eqref{Hua}:
\begin{equation*}
R_M(a_1,\dots,a_n;q)=
\sum_{\la^{(1)},\dots,\la^{(n)}}
\frac{q^{\frac{1}{2}\sum_{i,j=1}^n C_{ij}\pp{{\la^{(i)}}'}{{\la^{(j)}}'}}
\prod_{i=1}^n a_i^{\abs{\la^{(i)}}}}
{\prod_{i=1}^n (q;q)_{M_i-\ell(\la^{(i)})}b_{\la^{(i)}}(q)}.
\end{equation*}

By construction $R_M(a_1,\dots,a_n;q)$
satisfies the following invariance property.
\begin{lemma}\label{lemmainv}
We have
\begin{equation*}
\sum_{r_1=0}^{M_1}\cdots \sum_{r_n=0}^{M_n} 
\frac{q^{\frac{1}{2}\sum_{i,j=1}^n C_{ij}r_ir_j}
\prod_{i=1}^n a_i^{r_i}}
{\prod_{i=1}^n(q;q)_{M_i-r_i}} R_r(a_1,\dots,a_n;q)=
R_M(a_1,\dots,a_n;q).
\end{equation*}
\end{lemma}

\begin{proof}
Take the definition of $R_M$ given above and
replace each of $\lambda^{(1)},\dots,\lambda^{(n)}$
by its conjugate.
Then introduce the non-negative integer $r_i$ and the
partition $\mu^{(i)}$ with largest part not exceeding $r_i$
through $\lambda^{(i)}=(r_i,\mu^{(i)}_1,\mu^{(i)}_2,\dots)$.
Since $b_{\la'}(q)=(q;q)_{r-\mu_1} b_{\mu'}(q)$ for
$\la=(r,\mu_1,\mu_2,\dots)$ this implies
the identity of the lemma after again replacing
each of $\mu^{(1)},\dots,\mu^{(n)}$
by its conjugate.
\end{proof}

Next is the observation that the left-hand side
of \eqref{PPs} corresponds to $R_{(n,m)}(a,b;q)$. Hence we
may reformulate the A$_2$ instance of Lemma~\ref{lemmainv}.
\begin{theorem}
For $M_1$ and $M_2$ non-negative integers
\begin{multline}\label{BL}
\sum_{r_1=0}^{M_1}\sum_{r_2=0}^{M_2} 
\frac{a^{r_1}b^{r_2}q^{r_1^2-r_1r_2+r_2^2}}
{(q;q)_{M_1-r_1}(q;q)_{M_2-r_2}} 
\frac{(abq;q)_{r_1+r_2}}{(q,aq,abq;q)_{r_1}(q,bq,abq;q)_{r_2}} \\
=\frac{(abq;q)_{M_1+M_2}}{(q,aq,abq;q)_{M_1}(q,bq,abq;q)_{M_2}}.
\end{multline}
\end{theorem}

To see how this leads to the A$_2$ Rogers--Ramanujan identity \eqref{RRA2}
and its higher moduli generalizations, let $k_1,k_2,k_3$ be integers
such that $k_1+k_2+k_3=0$.
Making the substitutions 
\begin{align*}
r_1&\to r_1-k_1-k_2,& a&\to aq^{k_2-k_3}, & M_1&\to M_1-k_1-k_2, \\ 
r_2&\to r_2-k_1,    & b&\to bq^{k_1-k_2}, & M_2&\to M_2-k_1,
\end{align*}
in \eqref{BL}, we obtain
\begin{multline*}
\sum_{r_1=0}^{M_1}\sum_{r_2=0}^{M_2}
\frac{a^{r_1}b^{r_2}q^{r_1^2-r_1r_2+r_2^2}}
{(q;q)_{M_1-r_1}(q;q)_{M_2-r_2}}  \\ \times
\frac{(abq)_{r_1+r_2}}{(q;q)_{r_1+k_3}(aq;q)_{r_1+k_2}(abq;q)_{r_1+k_1}
(q;q)_{r_2-k_1}(bq;q)_{r_2-k_2}(abq;q)_{r_2-k_3}} \\
=\frac{a^{k_1+k_2}b^{k_1}q^{\frac{1}{2}(k_1^2+k_2^2+k_3^2)}
(abq)_{M_1+M_2}}{(q;q)_{M_1+k_3}(aq;q)_{M_1+k_2}(abq;q)_{M_1+k_1}
(q;q)_{M_2-k_1}(bq;q)_{M_2-k_2}(abq;q)_{M_2-k_3}},
\end{multline*}
which is equivalent to the type-II A$_2$ Bailey lemma 
of \cite[Theorem 4.3]{ASW99}. Taking $a=b=1$ this simplifies to
\begin{multline}\label{drie}
\sum_{r_1=0}^{M_1}\sum_{r_2=0}^{M_2}
\frac{q^{r_1^2-r_1r_2+r_2^2}}
{(q;q)_{M_1-r_1}(q;q)_{M_2-r_2}(q;q)_{r_1+r_2}^2}  
\qbin{r_1+r_2}{r_1+k_1}\qbin{r_1+r_2}{r_1+k_2}\qbin{r_1+r_2}{r_1+k_3}\\
=\frac{q^{\frac{1}{2}(k_1^2+k_2^2+k_3^2)}}{(q)_{M_1+M_2}^2}
\qbin{M_1+M_2}{M_1+k_1}\qbin{M_1+M_2}{M_1+k_2}\qbin{M_1+M_2}{M_1+k_3}.
\end{multline}

The idea is now to apply this transformation to the A$_2$ Euler identity
\cite[Equation (5.15)]{ASW99}
\begin{multline}\label{EulerA2}
\sum_{k_1+k_2+k_3=0} q^{\frac{3}{2}(k_1^2+k_2^2+k_3^2)} \\
\times
\sum_{w\in S_3}\epsilon(w)\prod_{i=1}^3
q^{\frac{1}{2}(3k_i-w_i+i)^2-w_i k_i} 
\qbin{M_1+M_2}{M_1+3k_i-w_i+i}
=\qbin{M_1+M_2}{M_1},
\end{multline}
where $w\in S_3$ is a permutation of $(1,2,3)$ and
$\epsilon(w)$ denotes the signature of $w$.

Replacing $M_1,M_2$ by $r_1,r_2$ in \eqref{EulerA2}, 
then multiplying both sides by
\begin{equation*}
\frac{q^{r_1^2-r_1r_2+r_2^2}}
{(q;q)_{M_1-r_1}(q;q)_{M_2-r_2}(q;q)_{r_1+r_2}^2},
\end{equation*}
and finally summing over $r_1$ and $r_2$ using \eqref{drie}
(with $k_i\to 3k_i-w_i+i$), yields
\begin{multline}\label{it1}
\sum_{k_1+k_2+k_3=0} q^{\frac{3}{2}(k_1^2+k_2^2+k_3^2)}
\sum_{w\in S_3}\epsilon(w)\prod_{i=1}^3
q^{(3k_i-w_i+i)^2-w_i k_i} 
\qbin{M_1+M_2}{M_1+3k_i-w_i+i} \\
=\sum_{r_1=0}^{M_1}\sum_{r_2=0}^{M_2}
\frac{q^{r_1^2-r_1r_2+r_2^2}(q;q)_{M_1+M_2}^2}
{(q;q)_{M_1-r_1}(q;q)_{M_2-r_2}(q;q)_{r_1}(q;q)_{r_2}(q;q)_{r_1+r_2}}.
\end{multline}
Letting $M_1$ and $M_2$ tend to infinity, and using
the Vandermonde determinant
\begin{equation*}
\sum_{w\in S_3}\epsilon(w)\prod_{i=1}^3 x_i^{i-w_i}
=\prod_{1\leq i<j\leq 3}(1-x_jx_i^{-1})
\end{equation*}
with $x_i\to q^{7k_i+2i}$, gives
\begin{multline*}
\frac{1}{(q;q)_{\infty}^3}\sum_{k_1+k_2+k_3=0} 
q^{\frac{21}{2}(k_1^2+k_2^2+k_3^2)-k_1-2k_2-3k_3} \\
\times
(1-q^{7(k_2-k_1)+2})(1-q^{7(k_3-k_2)+2})(1-q^{7(k_3-k_1)+4}) \\
=\sum_{r_1,r_2=0}^{\infty}
\frac{q^{r_1^2-r_1r_2+r_2^2}}
{(q;q)_{r_1}(q;q)_{r_2}(q;q)_{r_1+r_2}}.
\end{multline*}
Finally, by the A$_2$ Macdonald identity \cite{Macdonald72}
\begin{multline*}
\sum_{k_1+k_2+k_3=0}
\prod_{i=1}^3 x_i^{3k_i} q^{\frac{3}{2}k_i^2-ik_i}
\prod_{1\leq i<j\leq 3}(1-x_jx_i^{-1}q^{k_j-k_i}) \\
=(q;q)^2_{\infty}\prod_{1\leq i<j\leq 3}(x_i^{-1}x_j,q x_ix_j^{-1};q)_{\infty}
\end{multline*}
with $q\to q^7$ and $x_i\to q^{2i}$ this becomes
\begin{equation*}
\sum_{r_1,r_2=0}^{\infty}
\frac{q^{r_1^2-r_1r_2+r_2^2}}
{(q;q)_{r_1}(q;q)_{r_2}(q;q)_{r_1+r_2}}
=\frac{(q^2,q^2,q^3,q^4,q^5,q^5,q^7,q^7;q^7)_{\infty}}{(q;q)_{\infty}^3}.
\end{equation*}
This result is easily recognized as the A$_2$ Rogers--Ramanujan 
identity \eqref{RRA2}. 

The identity \eqref{it1} can be further iterated using \eqref{drie}.
Doing so and repeating the above calculations (requiring the
Vandermonde determinant with $x_i\to q^{(3n+1)k_i+ni}$
and the Macdonald identity with $q\to q^{3n+1}$ and $x_i\to q^{ni}$)
yields the following A$_2$ Rogers--Ramanujan-type identity for modulus 
$3n+1$ \cite[Theorem 5.1; $i=k$]{ASW99}:
\begin{multline*}
\sum_{\substack{\la,\mu \\ \ell(\la),\ell(\mu)\leq n-1}}
\frac{q^{\pp{\la}{\la}+\pp{\mu}{\mu}-\pp{\la}{\mu}}}
{b_{\la'}(q)b_{\mu'}(q)(q;q)_{\la_{n-1}+\mu_{n-1}}} \\
=\frac{(q^n,q^n,q^{n+1},q^{2n},q^{2n+1},
q^{2n+1},q^{3n+1},q^{3n+1};q^{3n+1})_{\infty}}{(q;q)_{\infty}^3}.
\end{multline*}
In the large $n$ limit ones recovers the A$_2$ case of Hua's
identity \eqref{Hua} with $a_1=a_2=1$.

To obtain identities corresponding to the modulus $3n-1$ we replace
$q\to 1/q$ in \eqref{EulerA2} to get 
\begin{multline*}
\sum_{k_1+k_2+k_3=0} q^{-\frac{3}{2}(k_1^2+k_2^2+k_3^2)}
\sum_{w\in S_3}\epsilon(w)\prod_{i=1}^3
q^{\frac{1}{2}(3k_i-w_i+i)^2+w_i k_i} 
\qbin{M_1+M_2}{M_1+3k_i-w_i+i} \\
=q^{2M_1M_2}\qbin{M_1+M_2}{M_1}.
\end{multline*}
Iterating this using \eqref{drie} and then taking the limit 
of large $M_1$ and $M_2$ yields \cite[Theorem 5.3; $i=k$]{ASW99}
\begin{multline*}
\sum_{\substack{\la,\mu \\ \ell(\la),\ell(\mu)\leq n-1}}
\frac{q^{\pp{\la}{\la}+\pp{\mu}{\mu}-\pp{\la}{\mu}+2\la_{n-1}\mu_{n-1}}}
{b_{\la'}(q)b_{\mu'}(q)(q;q)_{\la_{n-1}+\mu_{n-1}}} \\
=\frac{(q^{n-1},q^n,q^n,q^{2n-1},q^{2n-1},
q^{2n},q^{3n-1},q^{3n-1};q^{3n-1})_{\infty}}{(q;q)_{\infty}^3}.
\end{multline*}

Finally, the modulus $3n$ arises by iterating \cite[Equation (6.18)]{GK97}
\begin{equation*}
\sum_{k_1+k_2+k_3=0} \sum_{w\in S_3}\epsilon(w)\prod_{i=1}^3
q^{\frac{1}{2}(3k_i-w_i+i)^2} 
\qbin{M_1+M_2}{M_1+3k_i-w_i+i}
=\qbin{M_1+M_2}{M_1}_{q^3}.
\end{equation*}
A repeat of the earlier calculation then gives \cite[Theorem 5.4; $i=k$]{ASW99}
\begin{multline*}
\sum_{\substack{\la,\mu \\ \ell(\la),\ell(\mu)\leq n-1}}
\frac{q^{\pp{\la}{\la}+\pp{\mu}{\mu}-\pp{\la}{\mu}}
(q;q)_{\la_{n-1}}(q;q)_{\mu_{n-1}}}
{b_{\la'}(q)b_{\mu'}(q)(q;q)^2_{\la_{n-1}+\mu_{n-1}}}
\qbin{\la_{n-1}+\mu_{n-1}}{\la_{n-1}}_{q^3} \\
=\frac{(q^n,q^n,q^n,q^{2n},q^{2n},q^{2n},q^{3n},q^{3n};q^{3n})_{\infty}}
{(q;q)_{\infty}^3}.
\end{multline*}

\section{Some open problems}\label{secF}
In this final section we pose several open problems
related to the results of this paper.

\subsection{Macdonald's symmetric function}
The Hall--Littlewood functions $P_{\la}(x;t)$ 
are special cases of Macdonald's celebrated symmetric 
functions $P_{\la}(x;q,t)$, obtained from the latter by taking $q=0$.
An obvious question is whether Theorem~\ref{mainthm} can
be generalized to the Macdonald case.

From the Cauchy identity \cite[Sec. VI, Eqn. (4.13)]{Macdonald95}
\begin{equation*}
\sum_{\la}P_{\la}(x;q,t)Q_{\la}(y;q,t)=
\prod_{i,j\geq 1}\frac{(tx_iy_j;q)_{\infty}}{(x_iy_j;q)_{\infty}}
\end{equation*}
(see \cite{Macdonald95} for definitions related to Macdonald's symmetric function)
and the specialization
\begin{equation*}
Q_{\la}(1,t,t^2,\dots;q,t)=\frac{t^{n(\la)}}{c'_{\la}(q,t)}
\end{equation*}
we have
\begin{equation}\label{nPsom2}
\sum_{\la}\frac{t^{n(\la)}P_{\la}(x;q,t)}{c'_{\la}(q,t)}
=\prod_{i\geq 1}\frac{1}{(x_i;q)_{\infty}}.
\end{equation}
Here $c'_{\la}$ is the generalized hook-polynomial
\begin{equation*}
c'_{\la}(q,t)=\prod_{s\in\la}(1-q^{a(s)+1}t^{\ell(s)})
\end{equation*}
with $a(s)=\la_i-j$ and $\ell(s)=\la_j'-i$ the arm-length and leg-length
of the square $s=(i,j)$ of $\lambda$. Note that $c'_{\la}(0,t)=1$.

In view of the above we pose the problem of finding a $(q,t)$-analogue
of Theorem~\ref{mainthm} which simplifies to \eqref{nPsom2} when $y_i=0$
for all $i\geq 1$ and to \eqref{Pid} (with $q\to t$) when $q=0$.

Alternatively we may ask for a $(q,t)$-analogue of
\eqref{psiphi}. From \eqref{nPsom2} and standard properties of Macdonald
polynomials it follows that 
\begin{equation}\label{psiqt}
\sum_{\substack{\nu\supset\mu \\ \nu-\mu \textup{ hor. strip}}}
\frac{t^{n(\nu)} z^{\abs{\nu-\mu}}\psi_{\nu/\mu}(q,t)}{c'_{\nu}(q,t)}
=\frac{1}{(z;q)_{\infty}} \frac{t^{n(\mu)}}{c'_{\mu}(q,t)}.
\end{equation}
Here $\psi_{\la/\mu}(q,t)$ is generalization of $\psi_{\la/\mu}(t)$
(such that $\psi_{\la/\mu}(0,t)=\psi_{\la/\mu}(t)$)
given by
\begin{equation*}
\psi_{\la/\mu}(q,t)=\sideset{}{'}\prod_{s\in\mu} \frac{b_{\mu}(s)}{b_{\la}(s)},
\end{equation*}
where the product is over all squares $s=(i,j)$ of $\mu$ such that
$\theta_i>0$ and $\theta'_j=0$ for $\theta=\la-\mu$.
Moreover
\begin{equation*}
b_{\la}(s)=\frac{1-q^{a(s)}t^{l(s)+1}}{1-q^{a(s)+1}t^{l(s)}}
\end{equation*}
for $s\in\la$.

Hence a $(q,t)$-version of \eqref{psiphi} should reduce to \eqref{psiqt}
when $\la=0$ and to \eqref{psiphi} (with $q\to t$) when $q=0$.
Moreover, its right-hand side should involve the rational function
\begin{equation*}
\phi_{\la/\mu}(q,t)=\sideset{}{'}\prod_{s\in\la} \frac{b_{\mu}(s)}{b_{\la}(s)},
\end{equation*}
where the product is over all squares $s=(i,j)$ of $\la$ such that
$\theta_i'>0$ with $\theta=\la-\mu$ (and $b_{\mu}(s)=1$ if $s\not\in\mu$).

\subsection{The A$_n$ version of Theorem~\ref{mainthm}}
In the introduction we already mentioned the problem of evaluating the
A$_n$ sum \eqref{Ansum}. For $n>2$ this sum does not
equate to an infinite product and a possible scenario is that
for general $n$ the right-hand takes the form of an $n$ by $n$
determinant with infinite-product entries.

A specialized case of the sum \eqref{Ansum} does however exhibit a simple closed 
form evaluation, and the following extension of Theorem~\ref{mainthm} holds.
Let $x^{(1)}=x=(x_1,x_2,\dots)$, $x^{(n)}=y=(y_1,y_2,\dots)$
and $x^{(i)}_j=a_i q^j$ for $2\leq i\leq n-1$ and $j\geq 1$.
Also, let $\Delta_{+}'$ be the set (of cardinality $\binom{n-1}{2}$)
of positive roots of A$_n$ not containing the simple roots 
$\alpha_1$ and $\alpha_n$, i.e., the
set of roots of the form 
$\alpha_i+\alpha_{i+1}+\cdots+\alpha_j$ 
with $2\leq i\leq j\leq n-1$. Then
\begin{multline}
\sum_{\la^{(1)},\dots,\la^{(n)}} 
\prod_{i=1}^n q^{n(\la^{(i)})-\pp{{\la^{(i)}}'}{{\la^{(i+1)}}'}}
P_{\la^{(i)}}(x^{(i)};q) \\
=\prod_{\alpha\in\Delta_{+}'}
\frac{1}{(a^{\alpha}q;q)_{\infty}}
\prod_{i\geq 1}\prod_{j=1}^{n-1}
\frac{1}{(1-a_2\cdots a_j x_i)(1-a_{n-j+1}\cdots a_{n-1}y_i)}
\\ \times
\prod_{i,j\geq 1}\frac{1-a_2\cdots a_{n-1}x_iy_j}
{1-q^{-1}a_2\cdots a_{n-1}x_iy_j}.
\end{multline}
When $x_i=a_1q^i$ and $y_i=a_nq^i$ for $i\geq 1$ this yields
\eqref{Hua}.

Similarly, we have an isolated result for A$_3$ of the form
\begin{multline}
\sum_{\la,\mu,\nu} q^{n(\la)+n(\mu)+n(\nu)
-\pp{\la'}{\mu'}-\pp{\mu'}{\nu'}} \\
\times
P_{\la}(aq,aq^2,\dots;q)P_{\mu}(x;q)P_{\nu}(bq,bq^2,\dots;q) \\[2mm]
=\frac{1}{(aq,bq;q)_{\infty}}
\prod_{i\geq 1}\frac{(1-abx_i^2)}{(1-x_i)(1-a x_i)(1-bx_i)(1-abx_i)}
\prod_{i<j}\frac{1-abx_ix_j}{1-q^{-1}abx_ix_j}.
\end{multline}

\subsection{Bounds on Theorem~\ref{mainthm}}

The way we have applied \eqref{Pid} to obtain the A$_2$ Rogers--Ramanujan 
identity \eqref{RRA2} is rather different from Stembridge's Hall--Littlewood 
approach to the classical Rogers--Ramanujan identities \cite{Stembridge90}.
Specifically, Stembridge took \cite[p. 231]{Macdonald95}
\begin{equation*}
\sum_{\la} P_{2\la}(x;q)=\prod_{i=1}^n\frac{1}{1-x_i^2}\prod_{1\leq i<j\leq n}
\frac{1-qx_ix_j}{1-x_ix_j}=:\Psi(x;q)
\end{equation*}
for $x=(x_1,\dots,x_n)$, and generalized this to
\begin{equation}\label{stem}
\sum_{k=0}^{\infty}u^k \sum_{\substack{\la \\ \la_1\leq k}} 
P_{2\la}(x;q)=\sum_{\epsilon\in\{-1,1\}^n}
\frac{\Psi(x^{\epsilon};q)}{1-u x^{1-\epsilon}},
\end{equation}
where $f(x^{\epsilon})=f(x_1^{\epsilon_1},\dots,x_n^{\epsilon_n})$
and $x^{1-\epsilon}=x_1^{1-\epsilon_1} \cdots x_n^{1-\epsilon_n}$.
By specializing $x_i=z^{1/2}q^{i-1}$ for all $1\leq i\leq n$ this yields
\begin{multline}\label{St}
\sum_{\substack{\la \\ \la_1\leq k}}
\frac{z^{\abs{\la}}q^{2n(\la)}(q;q)_n}{(q;q)_{n-\ell(\la)}b_{\lambda}(q)} \\
=\sum_{r=0}^n (-1)^r (1-zq^{2r-1})
z^{(k+1)r} q^{(2k+3)\binom{r}{2}} \qbin{n}{r}
\frac{(z/q;q)_r}{(z/q;q)_{n+r+1}}.
\end{multline}
Letting $n$ tend to infinity and taking $k=1$ and $z=q$ or $z=q^2$ 
gives the Rogers--Ramanujan identities \eqref{RR} by an appeal
to the Jacobi triple-product identity to transform the sum on the right
into a product.

An obvious question is whether the identity \eqref{Pid} also admits
a version in which the partitions $\la$ and $\mu$ are summed
restricted to $\la_1\leq k_1$ and $\mu\leq k_2$, and if so, whether
such an identity would yield further A$_2$ $q$-series identities upon
specialization.
At present we have been unable to answer these questions. It is to be
noted, however, that since \eqref{nPsum} is the special case of \eqref{Pid} 
--- obtained by setting all $y_i$ equal to zero --- 
a bounded form of \eqref{nPsum} would be a precursor
to a bounded form of \eqref{Pid}.

Defining
\begin{equation*}
\Phi(x;q)=\prod_{i=1}^n\frac{1}{1-x_i}
\prod_{1\leq i<j\leq n}\frac{1-qx_ix_j}{1-x_ix_j}
\end{equation*}
Macdonald has shown that \cite[p.231--234]{Macdonald95}
\begin{equation}\label{Pphi}
\sum_{\la}P_{\la}(x;q)=\Phi(x;q)
\end{equation}
and
\begin{equation}\label{Md}
\sum_{k=0}^{\infty} u^k \sum_{\substack{\la \\ \la_1\leq k}} P_{\la}(x;q)
=\sum_{\epsilon\in\{-1,1\}^n}
\frac{\Phi(x^{\epsilon};q)}{1-u x^{(1-\epsilon)/2}}.
\end{equation}
With the above notation, \eqref{nPsum} takes a form rather similar 
to \eqref{Pphi};
\begin{equation*}
\sum_{\la}q^{n(\la)}P_{\la}(x;q)=\Phi(x;1).
\end{equation*}
But more can be done as the following bounded analogue of
\eqref{nPsum} holds.
\begin{theorem}\label{thmpf}
Let $[n]=\{1,\dots,n\}$.
For $I$ a subset of $[n]$ let $\abs{I}$ be its cardinality and
$J=[n]-I$ its complement. Then
\begin{multline}\label{idthm}
\sum_{k=0}^{\infty}u^k \sum_{\substack{\la \\ \la_1\leq k}}
q^{n(\la)}P_{\la}(x;q)  \\
=\sum_{I\subset [n]}
\frac{1}{1-u q^{\binom{\abs{I}}{2}}\prod_{i\in I}x_i}
\prod_{i\in I}\frac{1}{1-x_i^{-1}q^{1-\abs{I}}}
\prod_{j\in J}\frac{1}{1-x_j q^{\abs{I}}}
\prod_{\substack{i\in I \\ j\in J}}\frac{x_i-q x_j}{x_i-x_j}.
\end{multline}
\end{theorem}
If we specialize $x_i=zq^{i-1}$ --- but do not yet use
\eqref{Pspec} --- and equate coefficients of $u^k$, this 
leads to 
\begin{multline}\label{St2}
\sum_{\substack{\la \\ \la_1\leq k}}
q^{n(\la)}z^{\abs{\la}}P_{\la}(1,q,\dots,q^{n-1};q) \\
=\sum_{r=0}^n (-1)^r (1-zq^{2r-1})
z^{(k+1)r} q^{(2k+3)\binom{r}{2}} \qbin{n}{r}
\frac{(z/q;q)_r}{(z/q;q)_{n+r+1}}.
\end{multline}
This is a finite-$n$ analogue of \cite[Theorem 2]{Fulman00} of Fulman.
(To get Fulman's theorem take $z=q$ or $q^2$, replace $k\to k-1$, $q\to q^{-1}$ and
let $n$ tend to infinity. The Jacobi triple-product identity does the rest).
However, the reader should also note that the above right-hand side
coincides with the right-hand side of Stembridge's \eqref{St}.
Indeed, using the specialization formula \eqref{Pspec}, \eqref{St2}
is readily seen to be equivalent to \eqref{St} ---
the reason for this coincidence being that
$P_{2\la}(z^{1/2},z^{1/2}q,\dots,z^{1/2}q^{n-1};q)=
q^{n(\la)}P_{\la}(z,zq,\dots,zq^{n-1};q)$.

\begin{proof}[Proof of \eqref{St2}]
The left-hand side simply follows by extracting
the coefficient of $u^k$ in \eqref{idthm} and by making
the required specialization.

To get to the claimed right-hand side we note that
after specialization the term
\begin{equation*}
\prod_{\substack{i\in I \\ j\in J}}\frac{x_i-q x_j}{x_i-x_j}
\end{equation*} 
will vanish if there is an $i\in I$ and a $j\in J$ such that
$i-j=1$. Hence the only $I$ that will
contribute to the sum are the sets $\{1,\dots,r\}$ with $0\leq r\leq n$,
resulting in
\begin{multline*}
\sum_{r=0}^n \qbin{n}{r}
\frac{1}{(z^{-1}q^{2-2r};q)_r (zq^{2r};q)_{n-r}}
\frac{1}{1-uz^r q^{2\binom{r}{2}}}  \\
=\sum_{r=0}^n (-1)^r (1-zq^{2r-1})
z^r q^{3\binom{r}{2}}\qbin{n}{r}
\frac{(z/q;q)_r}{(z/q;q)_{n+r+1}}
\frac{1}{1-uz^r q^{2\binom{r}{2}}}.
\end{multline*}
The observation that the coefficient of $u^k$ of this series is given by 
the right-hand side of \eqref{St2} completes the proof.
\end{proof}

\begin{proof}[Proof of Theorem~\ref{thmpf}]
The proof proceeds along the lines of Macdonald's 
partial fraction proof of \eqref{Md} \cite{Macdonald95} 
and Stembridge's proof of \eqref{stem} (see also \cite{IJZ04,JZ01,JZ04}.

For any subset $E$ of $X=\{x_1,\dots,x_n\}$, let $p(E)$ denote
the product of the elements of $E$.
Let $\la=(\la_1,\dots,\la_n)$ be of the form 
$(\mu_1^{r_1},\dots,\mu_k^{r_k})$, with $\mu_1>\mu_2>\dots>\mu_k\geq 0$
and $r_1,\dots,r_k>0$ such that $\sum_i \mu_i=n$.
Then the defining expression \eqref{Pdef} of the Hall--Littlewood polynomials
can be rewritten as
\begin{equation}\label{Pf}
P_{\la}(x;q)=\sum_f \prod_{i=1}^k p(f^{-1}(i))^{\mu_i} 
\prod_{f(x_i)<f(x_j)}\frac{x_i-qx_j}{x_i-x_j},
\end{equation}
where the sum is over all surjections $f:X\to\{1,\dots,k\}$ such that
$\abs{f^{-1}(i)}=r_i$.
Each such surjection $f$ corresponds to a filtration 
\begin{equation*}
\F: \emptyset=F_0\subset F_1\subset \cdots\subset F_k=X,
\end{equation*}
according to the rule that $x\in F_i$ iff $f(x)\leq i$, and each
filtration of length $k$ such that $\abs{F_i-F_{i-1}}=r_i$ corresponds to
a surjection $f:X\to\{1,\dots,k\}$ such that $\abs{f^{-1}(i)}=r_i$. 
Hence \eqref{Pf} can be put as
\begin{equation*}
P_{\la}(x;q)=\sum_{\F} \pi_{\F}(X)\prod_{i=1}^k p(F_i-F_{i-1})^{\mu_i}
=\sum_{\F} \pi_{\F}(X)\prod_{i=1}^k p(F_i)^{\mu_i-\mu_{i+1}}
\end{equation*}
Here $\mu_{k+1}:=0$, the sum over $\F$ is a sum over all filtrations
of length $k$ such that $\abs{F_i-F_{i-1}}=r_i$, and
\begin{equation*}
\pi_{\F}(X)=\prod_{f(x_i)<f(x_j)}\frac{x_i-qx_j}{x_i-x_j},
\end{equation*}
with $f$ the surjection corresponding to $\F$.

Now given $\la$, the statistic $n(\la)$ may be expressed in terms
of the above defined quantities as
\begin{equation*}
n(\la)=\sum_{i=1}^k (\mu_i-\mu_{i+1})\binom{\abs{F_i}}{2}.
\end{equation*}
Hence, denoting the sum on the left of \eqref{idthm} by $S(u)$,
\begin{equation*}
S(u)=\sum_{\F}\pi_{\F}(X) \sum u^k
\prod_{i=1}^k 
\Bigl(q^{\binom{\abs{F_i}}{2}} p(F_i)\Bigr)^{\mu_i-\mu_{i+1}},
\end{equation*}
where the sum over $\F$ is a sum over filtrations of arbitrary length $k$
and where the inner sum is a sum over integers $k',\mu_1,\dots,\mu_k$ such that
$k'\geq\mu_1>\cdots>\mu_k\geq 0$.
Introducing the new variables $\nu_0=k'-\mu_1$ and $\nu_i=\mu_i-\mu_{i+1}$
for $i\in\{1,\dots,k\}$, so that $\nu_0,\nu_k\geq 0$ and all other $\nu_i>0$,
the inner sum can readily be carried out yielding
\begin{equation}\label{SA}
S(u)=\frac{1}{1-u}\sum_{\F}\pi_{\F}(X) A_{\F}(X,u),
\end{equation}
with
\begin{equation}\label{A}
A_{\F}(X,u)=\frac{1}{1-p(X)q^{\binom{n}{2}}u}
\prod_{i=1}^{k-1} \frac{p(F_i)q^{\binom{\abs{F_i}}{2}}u}
{1-p(F_i)q^{\binom{\abs{F_i}}{2}}u}.
\end{equation}

In the remainder it will be convenient not to work 
with the filtrations $\F$ but with the filtrations $\G$
\begin{equation*}
\G: \emptyset=G_0\subset G_1\subset \cdots\subset G_k=[n],
\end{equation*}
where $\G$ is determined from $\F$ by
$G_i=\{j|x_j\in F_i\}$.
Instead of $\pi_{\F}$ and $A_{\F}$ we will write
$\pi_{\G}$ and $A_{\G}$ and so on.

From \eqref{SA} and \eqref{A} it follows that the
following partial fraction expansion for $S(u)$ must hold:
\begin{equation*}
S(u)=\sum_{I\subset[n]}\frac{a_{I}}
{1-x_I q^{\binom{\abs{I}}{2}}u},
\end{equation*}
where $x_I$ stands for $\prod_{i\in I} x_i$.
After comparing this with \eqref{idthm}, the remaining task is to 
show that 
\begin{align}\notag
a_I&=\lim_{u\to x_I^{-1}q^{-\binom{\abs{I}}{2}}}
(1-x_I q^{\binom{\abs{I}}{2}}u) S(u) \\
&=\lim_{u\to x_I^{-1}q^{-\binom{\abs{I}}{2}}}
\frac{1-x_I q^{\binom{\abs{I}}{2}}u}{1-u} 
\sum_{\G} \pi_{\G}(X) A_{\G}(X,u) \label{sumG}
\end{align}
is given by
\begin{equation}\label{aI}
a_I=\prod_{i\in I}\frac{1}{1-x_i^{-1}q^{1-\abs{I}}}
\prod_{j\in J}\frac{1}{1-x_j q^{\abs{I}}}
\prod_{\substack{i\in I \\ j\in J}}\frac{x_i-q x_j}{x_i-x_j}.
\end{equation}

Since
\begin{equation*}
S(u)=\sum_{\la} q^{n(\la)}P_{\la}(x;q)
\sum_{k=\la_1}^{\infty}u^k=
\frac{1}{1-u}\sum_{\la} u^{\la_1} q^{n(\la)}P_{\la}(x;q),
\end{equation*}
we have
\begin{equation}\label{a0}
a_{\emptyset}=\lim_{u\to 1} (1-u) S(u)
=\sum_{\G} \pi_{\G}(X) A_{\G}(X,1)=\Phi(X)=\Phi(x;1),
\end{equation}
where, for later reference, we have introduced
\begin{equation*}
\Phi(Y)=\prod_{y\in Y}\frac{1}{1-y}
\end{equation*}
for arbitrary sets $Y$.
Now let us use \eqref{a0} to compute $a_I$ for general sets $I$.

The only filtrations that contribute to the sum in
\eqref{sumG} are those $\G$ that contain a $G_r$ (with
$0\leq r\leq k$) such that $G_r=I$.
Any such $\G$ may be decomposed into two filtrations 
$\G_1$ and $\G_2$ of length $r$ and $k-r$ by 
\begin{equation*}
\G_1:\emptyset=I-G_r\subset I-G_{r-1}\subset \cdots \subset I-G_1\subset I-G_0=I
\end{equation*}
and
\begin{equation*}
\G_2:\emptyset=G_r-I \subset G_{r+1}-I\subset \cdots \subset G_{k-1}-I
\subset G_k-I=[n]-I=J,
\end{equation*}
and given $\G_1$ and $\G_2$ we can clearly reconstruct $\G$.

For fixed $I$ and $J=[n]-I$ let $\bar{X}_I,X_J\subset X$ be the sets 
$\{x_i^{-1}q^{1-\abs{I}}|i\in I\}$ and
$\{x_jq^{\abs{I}}|j\in J\}$, respectively.
Then it is not hard to verify that
\begin{equation*}
\pi_{\G}(X)=\pi_{\G_1}(\bar{X}_I)\pi_{\G_2}(X_J)
\prod_{\substack{i\in I \\ j\in J}}\frac{x_i-qx_j}{x_i-x_j}.
\end{equation*}
Here we should perhaps remark that due to the homogeneity of the
terms making up $\pi_{\G}(X)$, the factors
$q^{1-\abs{I}}$ and $q^{\abs{I}}$ occurring in the definitions of 
$\bar{X}_I$ and $X_J$ simply cancel out.
Similarly, it follows that
\begin{equation*}
\frac{1-x_I q^{\binom{\abs{I}}{2}}u}{1-u}
A_{\G}(X,u)=x_Iq^{\binom{\abs{I}}{2}}u
\times A_{\G_1}(\bar{X}_I,x_I q^{\binom{\abs{I}}{2}}u)
\times A_{\G_2}(X_J,x_I q^{\binom{\abs{I}}{2}}u).
\end{equation*}

Substituting the above two decompositions in \eqref{sumG}
and taking the limit yields
\begin{align*}
a_I
&=\sum_{\G_1}  \pi_{\G_1}(\bar{X}_I) A_{\G_1}(\bar{X}_I,1)
\sum_{\G_2}  \pi_{\G_2}(X_J) A_{\G_2}(X_J,1)
\prod_{\substack{i\in I \\ j\in J}}\frac{x_i-qx_j}{x_i-x_j} \\
&=\Phi(\bar{X}_I)\Phi(X_J)
\prod_{\substack{i\in I \\ j\in J}}\frac{x_i-qx_j}{x_i-x_j} \\
&=\prod_{i\in I}\frac{1}{1-x_i^{-1}q^{1-\abs{I}}}
\prod_{j\in I}\frac{1}{1-x_j^{-1}q^{\abs{J}}}
\prod_{\substack{i\in I \\ j\in J}}\frac{x_i-qx_j}{x_i-x_j}
\end{align*}
in accordance with \eqref{aI}.
\end{proof}
Unfortunately, Macdonald's the partial fraction method fails to
provide an expression for
\begin{multline*}
\sum_{k_1,k_2=0}^{\infty}u^{k_1}v^{k_2}
\sum_{\substack{\la,\mu \\ \la_1\leq k_1,~\mu_1\leq k_2}} 
q^{n(\la)+n(\mu)-\pp{\la'}{\mu'}} P_{\la}(x;q)P_{\mu}(y;q)
\end{multline*}
when not all $y_i$ (or $x_i$) are equal to zero.

In fact, even for the special $1$-dimensional subcase
of Corollary~\ref{Cor3} no simple closed form expression 
is apparent for
\begin{equation*}
\sum_{k=0}^{\infty} u^k \sum_{\substack{\la \\ \la_1\leq k}} 
q^{n(\la)-\sum_{l=1}^j \la'_l} P_{\la}(x;q).
\end{equation*}

\bibliographystyle{amsplain}

\end{document}